\documentclass[12pt]{amsart}

\input xy
\xyoption{all}

\usepackage{cite}

\usepackage{geometry}
\usepackage{setspace}
\usepackage{amsmath}
\usepackage{amssymb}
\usepackage{amsthm}  
\newtheorem{same}{This should never appear}[section]
\newtheorem{defin}[same]{Definition}

\newtheorem{remark}[same]{Remark}
\newtheorem{theorem}[same]{Theorem}

\newtheorem{lemma}[same]{Lemma}
\newtheorem{fact}[same]{Fact}
\newtheorem{question}[same]{Question}

\newtheorem{prop}[same]{Proposition}

\newbox\noforkbox \newdimen\forklinewidth
\setbox0\hbox{$\textstyle\smile$}
\setbox1\hbox to \wd0{\hfil\vrule width \forklinewidth depth-2pt
 height 10pt \hfil}
\setbox\noforkbox\hbox{\lower 2pt\box1\lower 2pt\box0\relax}
\def\unionstick{\mathop{\copy\noforkbox}\limits}

\def\nonfork_#1{\unionstick_{\textstyle #1}}

\setbox0\hbox{$\textstyle\smile$}
\setbox1\hbox to \wd0{\hfil{\sl /\/}\hfil}
\setbox2\hbox to \wd0{\hfil\vrule height 10pt depth -2pt width
              \forklinewidth\hfil}
\newbox\doesforkbox
\setbox\doesforkbox\hbox{\lower 2pt\box1 \lower 2pt\box2\lower2pt\box0\relax}
\def\nunionstick{\mathop{\copy\doesforkbox}\limits}

\def\fork_#1{\nunionstick_{\textstyle #1}}

\newcommand{\sea}{\mathfrak{C}}

\newcommand{\ba}{\bold{a}}

\newcommand{\bx}{\bold{x}}
\newcommand{\by}{\bold{y}}

\newcommand{\dom}{\textrm{dom }}

\newcommand{\cf}{\text{cf }}
\newcommand{\rest}{\upharpoonright}

\newcommand{\bigM}{\widehat{M}}
\newcommand{\bigN}{\widehat{N}}
\newcommand{\bigf}{\widehat{f}}

\newcommand{\seq}[1]{\langle #1 \rangle}
\newcommand{\te}[1]{\text{#1}}

\title{Computing the number of types of infinite length}
\author{Will Boney}
\email{wboney@cmu.edu}
\address{Department of Mathematical Sciences \\ Carnegie Mellon University \\ Pittsburgh, Pennsylvania, USA}

\date{May 8, 2014} 

\begin{document}

\maketitle

\begin{abstract}
We show that the number of types of sequences of tuples of a fixed length can be calculated from the number of 1-types and the length of the sequences.  Specifically, if $\kappa \leq \lambda$, then
$$\sup_{\|M\| = \lambda} |S^\kappa(M)| = \left( \sup_{\|M\| = \lambda} |S^1(M)| \right)^\kappa$$
We show that this holds for any abstract elementary class with $\lambda$ amalgamation.  Basic examples show that no such calculation is possible for nonalgebraic types.  However, we introduce a generalization of nonalgebraic types for which the same upper bound holds.
\end{abstract}

\tableofcontents

\section{Introduction}

A well-known result in stability theory is that stability for 1-types implies stability for $n$-types for all $n < \omega$; see Shelah \cite{shelahfobook} Corollary I.\S2.2 or Pillay \cite{pillay}.0.9.  In this paper, we generalize this result to types of infinite length.

\begin{theorem} \label{introtheorem}
Given a complete theory $T$, if the supremum of the number of 1-types over models of size $\lambda \geq |T|$ is $\mu$, then for any (possibly finite) cardinal $\kappa \leq \lambda$, the supremum of the number of $\kappa$-types over models of size $\lambda$ is exactly $\mu^\kappa$.
\end{theorem}

We do this  by using the semantic, rather than syntactic, properties of types.  This allows our arguments to work in many nonelementary classes.  Thus, we work in the framework of Abstract Elementary Classes (AECs), which was introduced by Shelah in \cite{sh88}.  As we discuss in Section \ref{prelimsec}, AECs and Galois types include elementary classes and syntactic types and various nonelementary classes, such as those axiomatized in $L_{\lambda^+, \omega}(Q)$.  We use our results to answer a question of Shelah from \cite{shelahfobook}.

While the number of types of sequences of infinite lengths has not been calculated before, these types have already seen extensive use under the name $tp_*$ in \cite{shelahfobook} and $TP_*$ in \cite{shelahaecbook}.V.D.\S 3.  While \cite{shelahfobook} uses them most extensively, it is the use in \cite{shelahaecbook}.V.D.\S3 as types of models that might be most useful.  This means that stability in $\lambda$ can control the number of extensions of a model of size $\lambda$; see Section \ref{modelextension}.

After seeing preliminary versions of this work, Rami Grossberg asked if the above theorem could be proved for nonalgebraic types.  The examples in Proposition \ref{badnaex} show that such a theorem is not possible, even in natural elementary classes.  However, we introduce a generalization of nonalgebraic types of tuples called strongly separative types for which we can prove the same upper bound.  In AECs with disjoint amalgamation, such as elementary classes, nonalgebraic and strongly separative types coincide for types of length 1.  For longer types, we require that realizations are, in a sense, nonalgebraic over each other.  For instance, in $ACF_0$, the type of $(e, \pi)$ can be considered ``more nonalgebraic'' over the set of algebraic numbers than the type of $(e, 2e)$.  This is made precise in Definition \ref{sepdef}.

Finally, in Section \ref{satsec}, we investigate the saturation of types of various lengths.  The ``saturation = model homogeneity'' lemma (recall Lemma \ref{modhomsat}) shows that saturation is equivalent for all lengths.  We also use bounds on the number of types and various structural properties to construct saturated models.

This paper was written while working on a Ph.D. under the direction of Rami Grossberg at Carnegie Mellon University and I would like to thank Professor Grossberg for his guidance and assistance in my research in general and relating to this work specifically.  A preliminary version of this paper was presented in the CMU Model Theory Seminar and I'd like to thank the participants for helping to improve the presentation of the material, especially Alexei Kolesnikov.  He pointed out a gap in the proof of Theorem \ref{lowerbound} and vastly improved the presentation of Theorem \ref{finalone}, among many other improvements.  I would also like to thank my wife Emily Boney for her support.

\section{Preliminaries} \label{prelimsec}

We use the framework of abstract elementary classess to prove our results.  Thus, we offer the following primer on abstract elementary classes and Galois types.

The definition for an abstract elementary class (AEC) was first given by Shelah in \cite{sh88}.  The definitions and concepts in the section are all part of the literature; in particular, see Baldwin\cite{baldwinbook}, Shelah \cite{shelahaecbook}, Grossberg \cite{grossberg2002}, or the forthcoming Grossberg \cite{ramibook} for more information.

\begin{defin}
We say that $(K, \prec_K)$ is an Abstract Elementary Class iff
\begin{enumerate}

    \item There is some language $L = L(K)$ so that every element of $K$ is an $L(K)$-structure; 

    \item $\prec_K$ is a partial order on $K$;

    \item for every $M, N \in K$, if $M \prec_K N$, then $M \subseteq_{L} N$;

    \item $(K, \prec_K)$ respects $L(K)$ isomorphisms, if $f: N \to N'$ is an $L(K)$ isomorphism and $N \in K$, then $N' \in K$ and if we also have $M \in K$ with $M \prec_K N$, then $f(M) \in K$ and $f(M) \prec_K N'$;

    \item \emph{(Coherence)} if $M_0, M_1, M_2 \in K$ with $M_0 \prec_K M_2$, $M_1 \prec_K M_2$, and $M_0 \subseteq M_1$, then $M_0 \prec M_1$;

    \item \emph{(Tarski-Vaught axioms)} suppose $\seq{M_i \in K : i < \alpha}$ is a $\prec_K$-increasing continuous chain, then

        \begin{enumerate}

            \item $\cup_{i < \alpha} M_i \in K$ and, for all $i < \alpha$, we have $M_i \prec_K \cup_{i < \alpha} M_i$; and

            \item if there is some $N \in K$ so that, for all $i < \alpha$, we have $M_i \prec_K N$, then we also have $\cup_{i < \alpha} M_i \prec_K N$; and

        \end{enumerate}

    \item \emph{(L\"{o}wenheim-Skolem number)} $LS(K)$ is the minimal infinite cardinal $\lambda \geq |L(K)|$ such that for any $M \in K$ and $A \subset |M|$, there is some $N \prec_K M$ such that $A \subset |N|$ and $\|N|| \leq |A| + \lambda$.

\end{enumerate}

\end{defin}

\begin{remark}
We drop the subscript on $\prec_K$ when it is clear from context and we abuse notation by calling $K$ an AEC when we mean that $(K, \prec_K)$ is an AEC.  We follow the convention of Shelah and use $\| M \|$ to denote the cardinality of the universe of $M$.  In this paper, $K$ is always an AEC that has no models of size smaller than the L\"{o}wenheim-Skolem number.
\end{remark}

The most basic example of an AEC is any elementary class with the elementary substructure relation.  In particular, they are structural properties that can be proved without mention of the compactness theorem.  This means that classes of models axiomatized in many other logics, such as those with infinite conjunction/disjunction or with additional quantifiers, are also AECs with the appropriate substructure relation.  \cite{baldwinbook}.4 and \cite{genquant} give more details and \cite{nperp} discusses AECs consisting of modules.

We will briefly summarize some of the basic notations, definitions, and results for AECs.

\begin{defin} Let $K$ be an Abstract Elementary Class.
\begin{enumerate}

	\item \begin{eqnarray*}\vspace{-1in}
	K_\lambda &=& \{ M \in K : \|M\| = \lambda \}\\
	K_{\leq \lambda} &=& \{ M \in K : \|M\| \leq \lambda \}
	\end{eqnarray*}

	\item A $K$-embedding is an injection $f: M \to N$ that respects $L(K)$ such that $f(M) \prec N$.

	\item $K$ has the \emph{$\lambda$-amalgamation property} ($\lambda$-AP) iff for any $M \prec N_0, N_1 \in K_\lambda$, there is some $N^* \in K$ and $f_i: M \to N_i$ so that
\[
 \xymatrix{\ar @{} [dr] N_1  \ar[r]^{f_1}  & N^*\\
M \ar[u] \ar[r] & N_0 \ar[u]_{f_2}
 }
\]
commutes.

	\item $K$ has the \emph{$\lambda$-joint mapping property} ($\lambda$-JMP) iff for any $M_0, M_1 \in K_\lambda$, there is some $N \in K$ and $f_\ell:M_\ell \to N$ for $\ell = 0, 1$.
		
\end{enumerate}
\end{defin}

In cases where the AEC is axiomatized by a logic, the usefulness of types as sets of formulas comes from the unique features of first order logic such as compactness.  In order to compensate for this, Shelah isolated a semantic notion of type in \cite{sh300} that Grossberg named \emph{Galois type} in \cite{grossberg2002} that can replace sets of formulas.

We differ from the standard treatment of types in that we allow the length of our types to be possibly infinite.  This is necessary because we want to consider types of infinite tuples.

\begin{defin} \label{galoistype} Let $K$ be an AEC, $\lambda \geq LS(K)$, and $(I, <_I)$ an ordered set.
 \begin{enumerate}
  
    \item Set $K^{3, I}_\lambda = \{ (\seq{a_i : i \in I}, M, N) : M \in K_\lambda, M \prec N \in K_{\lambda + |I|}, \te{ and } \{a_i : i \in I \} \subset |N| \}$.  The elements of this set are referred to as \emph{pretypes}.

    \item Given two pretypes $(\seq{a_i : i \in I}, M, N)$ and $(\seq{b_i : i \in I}, M', N')$ from $K_\lambda^{3, I}$, we say that $(\seq{a_i : i \in I}, M, N) \sim_{AT} (\seq{b_i : i \in I}, M', N')$ iff $M = M'$ and there is $N^* \in K$ and $f: N \to N^*$ and $g: N' \to N^*$ so that $f(a_i) = g(b_i)$ for all $i \in I$ and the following diagram commutes:
\[
 \xymatrix{\ar @{} [dr] N'  \ar[r]^{g}  & N^*\\
M \ar[u] \ar[r] & N \ar[u]_{f}
 }
\]

    \item Let $\sim$ be the transitive closure of $\sim_{AT}$.
    
    \item For $M \in K$, set $gtp(\seq{a_i : i \in I}/M, N) = [(\seq{a_i : i \in I}, M, N)]_\sim$ and $gS^I(M) = \{ gtp(\seq{a_i : i \in I}/M, N) : (\seq{a_i : i \in I}/M, N) \in K^{3, I}_{\| M \|} \}$.  
    
    \item For $M \in K$, define $gS^I_{\te{na}}(M) = \{ tp(\seq{a_i : i \in I}/M, N) \in S^I(M) : a_i \in N - M \te{ for all $i \in I$} \}$.
    
    \item Let $M \in K$ and $p = gtp(\seq{a_i : i \in I}/M, N) \in gS^I(M)$.
    \begin{itemize}
    
    	\item If $M' \prec M$, then $p \rest M'$ is $gtp(\seq{a_i : i \in I}/M', N')$ for some (any) $N' \in K_{\|M'\|+ |I|}$ with $M' \prec N' \prec N$ and $\seq{a_i : i \in I} \subset |N'|$.
     
     	\item If $I_0 \subset I$, then $p^{I_0}$ is $gtp(\seq{a_i : i \in I_0}/M, N')$ for some (any) $N' \in K_{\|M\| + |I_0|}$ with $M \prec N' \prec N$ and $\seq{a_i : i \in I_0} \subset |N'|$.
	
	\end{itemize}

 \end{enumerate}

\end{defin}

\begin{remark}
If $K$ has the $\lambda + |I|$-amalgamation property, then $\sim_{AT}$ is transitive and, thus, an equivalence relation on $K^{3,I}_\lambda$.  Note that `$AT$' stands for ``atomic.''
\end{remark}

Since we will make extensive use of Galois types, we will assume that all AECs have the amalgamation property.  We will also use the joint mapping property as a ``connectedness'' property.  For first order theories, these properties follow from compactness and interpolation.

In the first order case, amalgamation over models follows directly from compactness and interpolation.  For complete theories, amalgamation holds over sets as well.  Furthermore, Galois types and syntactic types correspond.
This means that Theorem \ref{introtheorem} from the introduction follows from Theorem \ref{finalone} below and we can translate the other results similarly.  However, there is no AEC version of $tp_\Delta$ for $\Delta \subsetneq Fml(L)$; we summarize what we do know at the end of the next section.

For AECs axiomatized in other logics, the correspondence is not so nice.  Baldwin and Kolesnikov \cite{untame} analyzed the Hart-Shelah examples from \cite{hash323} to show that two elements can have the same syntactic type but different Galois types, even in an $L_{\omega_1, \omega}$ axiomatized class.

There is a correspondence in the other direction.  If the syntactic type of two elements are different in a logic that the AEC can ``see,'' then the Galois types must be different as well.  For instance, suppose $\psi$ is a sentence in some fragment $L_A$ of $L_{\lambda^+, \omega}$.  Then, if $tp_{L_A}(a/M, N_1) \neq tp_{L_A}(b/M, N_2)$, then their Galois types must differ in the AEC $(Mod \te{ } \psi, \prec_{L_A})$.  This means that that classical many-types theorems for non-first order logic, such as those for $L_{\omega_1, \omega}$ in \cite{keislerw1} and for $L(Q)$ in \cite{Keisler:1970fk}, imply many Galois types.

We investigate the supremum of the number of types of a fixed length over all models of a fixed size.  To simplify this discussion, we introduce the following notation.

\begin{defin} \label{tbdef}
The \emph{type bound for $\lambda$ sized domains and $\kappa$ lengths} is denoted $\mathfrak{tb}^\kappa_\lambda = \sup_{M \in K_\lambda} |gS^\kappa(M)|$.
\end{defin}

Shelah has introduced the notation of $tp_*$  in \cite{shelahfobook}.III.1.1 and $TP_*$ in \cite{shelahaecbook}.V.D.3 to denote the types of infinite tuples, with $tp_*$ having a syntactic definition (sets of formulas) and $TP_*$ having a  semantic definition (Galois types).  Thus, $\mathfrak{tb}_\lambda^\kappa$ counts the maximum number of types of a fixed length $\kappa$ over models of a fixed size $\lambda$, allowing for the possibility that this maximum is not achieved.  These long types are also used fruitfully in Makkai and Shelah \cite{makkaishelah}, Grossberg and VanDieren \cite{tamenessone}, and Boney and Grossberg \cite{shorttamedep}.

Clearly, $\lambda$-stability is the same as the statement that $\mathfrak{tb}_\lambda^1 = \lambda$.  Also, we always have $\mathfrak{tb}_\lambda^1 \geq \lambda$ because each element in a model has a distinct type.  Other notations have been used to count the supremum of the number of types, although the lengths have been finite.  In \cite{sixtheories}, Keisler uses
$$f_T(\kappa) = \sup \{ |S^1(M, N)| : M, N \models T, M \prec N, \te{ and } \|M\| = \kappa\}$$
In \cite{shelahfobook}.II.4.4, Shelah uses, for $\Delta \subset L(T)$ and $m < \omega$,
$$K^m_\Delta(\lambda, T) : = \min \{ \mu : |A| \leq \lambda \te{ implies } |S^m_\Delta(A)| < \mu\} = \sup_{|A| = \lambda} (|S^m_\Delta(A)|^+ ) $$
The relationships between these follow easily from the definitions
\begin{eqnarray*}
f_T(\kappa) &=& \mathfrak{tb}_\lambda^1\\
K_{L(T)}^m(\lambda, T) &=& \sup_{\|M\| = \lambda} ( |S^m(M)|^+) = 
\begin{cases}
\mathfrak{tb}_\lambda^m & \text{if $K_{L(T)}^m(\lambda, T)$ is limit}\\
(\mathfrak{tb}_\lambda^m)^+ & \text{if $K_{L(T)}^m(\lambda, T)$ is successor}\\
\end{cases}\\
&=& \begin{cases}
\mathfrak{tb}_\lambda^m & \text{if $\mathfrak{tb}^\lambda_m$ is a strict supremum}\\
(\mathfrak{tb}_\lambda^m)^+ & \text{if the supremum in $\mathfrak{tb}^\lambda_m$ is achieved}\\
\end{cases}\\
\end{eqnarray*}

From this last equality, a basic question concerning $\mathfrak{tb}^\kappa_\lambda$ is if the supremum is strict or if there is a model that achieves the value.  Below we describe two basic cases when the supremum in $\mathfrak{tb}^\kappa_\lambda$ is achieved.

\begin{prop} \label{achievesup}
Suppose $K$ is an AEC with $\lambda$-AP and $\lambda$-JMP and $\kappa \leq \lambda$.  If $\cf \mathfrak{tb}^\kappa_\lambda \leq \lambda$ or if $I(K, \lambda) \leq \lambda$, then there is $M \in K_\lambda$ such that $|gS^\kappa(M)| = \mathfrak{tb}^\kappa_\lambda$.
\end{prop}

{\bf Proof:} The idea of this proof is to put the $\leq \lambda$ many $\lambda$ sized models together into a single $\lambda$ sized model that will witness the conclusion.  Pick $\seq{M_i^* \in K_\lambda : i < \chi}$ with $\chi \leq \lambda$ such that $\{ |gS^\kappa(M_i^*)| : i < \chi\}$  has supremum $\mathfrak{tb}^\kappa_\lambda$; in the first case, this can be done by the definition of supremum and, in the second case, this can be done because there are only $I(K, \lambda)$ many possible values for $|gS^\kappa(M)|$ when $M \in K_\lambda$.  Using amalgamation and joint mapping, we construct increasing and continuous $\seq{N_i \in K_\lambda : i < \chi}$ such that $M_i^*$ is embeddable into $N_{i+1}$.  Set $M = \cup_{i < \chi} N_i$.  Since $\chi \leq \lambda$, we have $M \in K_\lambda$; this fact was also crucial in our construction.  Since $M_i^*$ can be embedded into $M$, we have that $|gS^\kappa(M_i^*)| \leq |gS^\kappa(M)| \leq \mathfrak{tb}^\kappa_\lambda$.  Taking the supremum over all $i < \chi$, we get $\mathfrak{tb}^\kappa_\lambda = |gS^\kappa(M)|$, as desired. \hfill \dag\\

The use of joint embedding here seems necessary, at least from a naive point of view.  It seems possible to have distinct AECs $K^n$ in a common language that have models $M^n \in K^n_\lambda$ such that $|gS^\kappa(M^n)| = \mathfrak{tb}^\kappa_\lambda = \lambda^{+n}$, each computed in $K^n$.  Then, we could form $K^\omega$ as the disjoint union of these classes; this would be an AEC with $\mathfrak{tb}^\kappa_\lambda = \lambda^{+\omega}$ and the supremum would not be achieved.  However, examples of such $K^n$, even with $\kappa=1$, are not known and the specified values of $|gS^\kappa(\cdot)|$ might not be possible.

These relationships help to shed light on a question of Shelah.

\begin{question}[\cite{shelahfobook}.III.7.6] \label{shelahqsol}
Is $K^m_{L(T)}(\lambda, T) = K^1_{L(T)}(\lambda, T)$ for $m < \omega$?
\end{question}

The answer is yes, even for a more general question, under some cardinal arithmetic assumptions.  Below, $\lambda^{(+\lambda^+)}$ denotes the $\lambda^+$th successor of $\lambda^+$.

\begin{theorem}
Suppose $2^\lambda < \lambda^{(+\lambda^+)}$.  If $\Delta \subset Fml(L(T))$ is such that $\phi(\bx, x, \by) \in \Delta$ implies $\exists z \phi(\bx, z, \by) \in \Delta$ and $n < \omega$, then
$$K^n_\Delta(\lambda, T) = K^1_\Delta(\lambda, T)$$
\end{theorem}

{\bf Proof:} There are two cases to consider: whether or not the supremum in $\mathfrak{tb}_\lambda^m$ is strict or is acheived.  If the supremum is strict, then we claim the supremum in $\mathfrak{tb}_\lambda^1$ is strict as well.  If not, there is some $M \in K_\lambda$ such that $|S^1(M)| = \mathfrak{tb}_\lambda^1$.  But then, by Theorem \ref{lowerbound},
$$\mathfrak{tb}_\lambda^m > |S^m(M)| \geq |S^1(M)|^m = (\mathfrak{tb}_\lambda^1)^m = \mathfrak{tb}_\lambda^m$$
a contradiction.  So $\mathfrak{tb}_\lambda^m$ is a strict supremum and 
$$K^m_{L(T)}(\lambda, T)= \mathfrak{tb}_\lambda^m = \mathfrak{tb}_\lambda^1  = K^1_{L(T)}(\lambda, T)$$
Note that this continues to hold if $m$ is infinite or if we consider the corresponding relationship for Galois types in an AEC with amalgamation.  Furthermore, this does not use the cardinal arithmetic assumption.

Now we consider the case that the supremum in $\mathfrak{tb}_\lambda^m$ is achieved and suppose for contradiction that the supremum in $\mathfrak{tb}_\lambda^1$ is strict.  Then $m > 1$ and we assume it is the minimal such $m$.  If $\mathfrak{tb}_\lambda^m = \mathfrak{tb}_\lambda^1$ is regular, than the pigeonhole argument used in Theorem \ref{lowerbound} can find a model achieving $\mathfrak{tb}_\lambda^1$.  In fact, this argument just requires that 
$$\sup \{ |S^{m-1}(Ma)| : a \vDash p, p \in S^1(M) \} < \lambda$$
By the remarks above the question, we know that $\cf \mathfrak{tb}_\lambda^1>\lambda$ since the supremum is strict.  This gives us that
$$\lambda < \cf \mathfrak{tb}_\lambda^1 < \mathfrak{tb}_\lambda^1 \leq 2^\lambda$$
However, this contradicts our cardinal arithmetic assumption because the minimal singular cardinal with cofinality above $\lambda$ is $\lambda^{(+\lambda^+)} > 2^\lambda$.  Thus
$$K^m_{L(T)}(\lambda, T)= (\mathfrak{tb}_\lambda^m)^+ = (\mathfrak{tb}_\lambda^1)^+  = K^1_{L(T)}(\lambda, T)$$
\hfill \dag\\

We now state the ``model-homogeneity = saturation'' lemma for AECs.  This has long been known for first order theories and first appeared for AECs in \cite{sh300}, although a correct proof was not given in print until Shelah \cite{sh576}.0.26.1.

\begin{lemma}[Shelah] \label{modhomsat}
Let $K$ be an AEC with amalgamation and $\lambda >  LS(K)$.  Then the following are equivalent for $M \in K$:
\begin{itemize}
	
	\item $M$ is $\lambda$-model homogeneous: for every $N_1 \prec N_2 \in K_{<\lambda}$ with $N_1 \prec M$, there is a $K$ embedding $f: N_2 \to_{N_1} M$; and
	
	\item $M$ is $\lambda$-Galois saturated: for every $N \prec M$ with $\|N\| < \lambda$ and every $p \in S^1(N)$, $p$ is realized in $M$.

\end{itemize} 
\end{lemma}

\section{Results on $S^\alpha(M)$}

This section aims to prove Theorem \ref{introtheorem} for AECs.  In our notation, this can be stated as follows.

\begin{theorem}\label{aectheorem}
If $K$ is an AEC with $\lambda$ amalgamation, then for any $\kappa \leq \lambda$, allowing $\kappa$ to be finite or infinite, we have $\mathfrak{tb}_\lambda^\kappa = (\mathfrak{tb}_\lambda^1)^\kappa$.
\end{theorem}

We prove this by proving a lower bound (Theorem \ref{lowerbound}) and an upper bound (Theorem \ref{finalone}) for $\mathfrak{tb}_\lambda^\kappa$.  Note that when $\kappa = \lambda$, this value is always the set-theoretic maximum, $2^\lambda$.  However, for $1 < \kappa < \min \{ \chi : (\mathfrak{tb}_\lambda^1)^\chi = 2^\lambda \}$, this provides new information.

For readers interested in AECs beyond elementary classes, we note the use of amalgamation for the rest of this section and for the rest of this paper.  It remains open whether these or other bounds can be found on the number of types without amalgamation.  One possible obstacle is that different types cannot be put together: if we assume amalgamation, then given two types $p, q \in gS^1(M)$, there is some type $r \in gS^2(M)$ such that its first coordinate extends $p$ and its second coordinate extends $q$.  This will be a crucial tool in the proof of the lower bound.  However, if we cannot amalgamate a model that realizes $p$ and a model that realizes $q$ over $M$, then such an extension type does not necessarily exist.

For the lower bound, we essentially ``put together'' all of the different types in $gS^1(M)$ as discussed above.

\begin{theorem} \label{lowerbound}
Let $K$ be an AEC with $\lambda$-AP and $\lambda$-JMP.  We have $\mathfrak{tb}_\lambda^\kappa \geq (\mathfrak{tb}_\lambda^1)^\kappa$.  In particular, given $M \in K_\lambda$, $|gS^\kappa(M)| \geq |gS^1(M)|^\kappa$.
\end{theorem}

{\bf Proof:}  We first prove the ``in particular'' clause and use that to prove the statement.  Fix $M \in K_\lambda$ and set $\mu = |gS^1(M)|$.  Fix some enumeration $\seq{p_i : i < \mu}$ of $gS^1(M)$.  Then we claim that there is some $M^+ \succ M$ that realizes all of the types in $gS^1(M)$.\\
To see this, let $N_i \succ M$ of size $\lambda$ contain a realization of $p_i$.  Then set $M_0 = M$ and $M_1 = N_0$.  For $\alpha = \beta + 1$, amalgamate $M_\beta$ and $N_\beta$ over $M$ to get $M_\alpha \succ M_\beta$ and $f: N_\beta \to_M M_\alpha$; since $N_\beta$ realizes $p_\beta \in S(M)$, $f(N_\beta)$ realizes $f(p_\beta) = p_\beta$.  So $M_\alpha$ does as well.  Take unions at limits.  Then $M^+ := \cup_{\beta < \alpha} M_\beta$ realizes each type in $gS^1(M)$.\\
Having proved the claim, we show that $|gS^\kappa(M)| \geq \mu^\kappa$.  For each $i < \mu$, pick $a_i \in |M^+|$ that realizes $p_i$.  For each $f \in {}^\kappa \mu$, set $\ba_f = \seq{ a_{f(i)} : i < \kappa}$.  We claim that the map $(f \in {}^\kappa \mu) \to gtp(\ba_f/M, M^+)$ is injective, which completes the proof.\\
To prove injectivity, note that $gtp(a_j/M, M^+) = gtp(a_k/M, M^+)$ iff $j = k$.  Suppose $gtp(\ba_f/M, M^+) = gtp(\ba_g/M, M^+)$.  Then, we see that $gtp(a_{f(i)}/M, M^+) = gtp(a_{g(i)}/M, M^+)$ for each $i < \kappa$.  By our above note, that means that $f(i) = g(i)$ for every $i \in \kappa = \dom f = \dom g$.  So $f = g$.  Thus, $|gS^\kappa(M)| \geq |{}^\kappa \mu| = \mu^\kappa$, as desired.  \\
Now we prove that $\mathfrak{tb}^\kappa_\lambda \geq (\mathfrak{tb}^1_\lambda)^\kappa$.  This is done by separating into cases based on $\cf (\mathfrak{tb}^1_\lambda)$.  If $\cf (\mathfrak{tb}^1_\lambda) > \kappa$, then it is known that exponentiating to $\kappa$ is continuous at $\mathfrak{tb}^1_\lambda$.  Stated more plainly, if $X$ is a set of cardinals such that $\cf (\sup_{\chi \in X} \chi) > \kappa$, then 
$$(\sup_{\chi \in X} \chi)^\kappa = \sup_{\chi \in X} (\chi^\kappa)$$
Then, we compute that
$$(\mathfrak{tb}^1_\lambda)^\kappa = (\sup_{M \in K_\lambda} |gS^1(M)|)^\kappa = \sup_{M \in K_\lambda} (|gS^1(M)|^\kappa) \leq \sup_{M \in K_\lambda} |gS^\kappa(M)| = \mathfrak{tb}^\kappa_\lambda$$
If $\cf( \mathfrak{tb}^1_\lambda )\leq \kappa$, then we also have $\cf (\mathfrak{tb}^1_\lambda )\leq \lambda$.  By Proposition \ref{achievesup}, we know that the supremum of $\mathfrak{tb}^1_\lambda$ is achieved, say by $M^* \in K_\lambda$.  Then
$$(\mathfrak{tb}^1_\lambda)^\kappa = |gS^1(M^*)|^\kappa \leq |gS^\kappa(M^*)| \leq \sup_{M \in K_\lambda} |gS^\kappa(M)| = \mathfrak{tb}^\kappa_\lambda$$
\hfill \dag \\

Now we show the upper bound.  We do this in two steps.  First, we present the ``successor step'' in Theorem \ref{firstone} to give the reader the flavor of the argument.  Then Theorem \ref{finalone} gives the full argument using direct limits.

\begin{theorem} \label{firstone}
For any AEC K with $\lambda$-AP and any $n < \omega$, $\mathfrak{tb}_\lambda^n \leq \mathfrak{tb}_\lambda^1$.
\end{theorem}

Note that, since it includes the $\| M \|$ many algebraic types, $gS^1(M)$ is always infinite, so this result could be written $\mathfrak{tb}_\lambda^n \leq (\mathfrak{tb}_\lambda^1)^n$.

{\bf Proof:} We prove this by induction on $n < \omega$.  The base case is $\mathfrak{tb}_\lambda^1 \leq \mathfrak{tb}_\lambda^1$.  Suppose $\mathfrak{tb}_\lambda^n \leq \mathfrak{tb}_\lambda^1$ and set $\mu = \mathfrak{tb}_\lambda^1$.  For contradiction, suppose there is some $M \in K_\lambda$ such that $|gS^{n+1}(M)| > \mu$.  Then we can find distinct $\{ p_i \in S^{n+1}(M) \mid i < \mu^+\}$ and find $\seq{a_j^i \mid j < n+1} \models p_i$ and $N_i \succ M$ that contains each $a^i_j$ for $j < n+1$.\\
Consider $\{ gtp(\seq{a_j^i \mid j < n} / M, N_i) : i < \mu^+ \} \subset gS^n(M)$.  By assumption, this set has size $\mu$.  So there is some $I \subset \mu^+$ of size $\mu^+$ such that, for all $i \in I$, $gtp(\seq{a_j^i \mid j < n}/M, N_i)$ is constant.\\
Fix $i_0 \in I$.  For any $i \in I$, the Galois types of $\seq{a_j^i : j < n}$ and $\seq{a_j^{i_0} : j < n}$ over $M$ are equal.  Thus, there are $N_i^* \succ N_{i_0}$ and $f_i: N_i \to_{M} N_i^* $ such that $f_i(a_j^i) = a_j^{i_0}$ for all $j < n$ and 
\[
\xymatrix{
N_{i_0} \ar[r] & N_i^* \\
M \ar[u] \ar[r] & N_i \ar[u]_{f_i}
}
\]
commutes.  Now consider the set $\{ gtp(f_i(a_n^i) / N_{i_0}, N_i^*) \mid i \in I\}$. We have that $|I| = \mu^+$ and $|S^1(N_{i_0})| \leq \mathfrak{tb}^1_\lambda = \mu$, so there is $I^* \subset I$ of size $\mu^+$ so, for all $i \in I^*$, $gtp(f_i(a_n^i) / N_{i_0}, N_i^*)$ is constant.  Let $i \neq k \in I^*$.\\
Then $gtp(f_i(a_n^i)/N_{i_0}, N_i^*) = gtp(f_k(a_n^k)/N_{i_0}, N_k^*)$.  By the definition of Galois types, we can find $N^{**}$, $g_k:N_k^* \to N^{**}$, and $g_i: N_i^* \to N^{**}$ such that $g_k(f_k(a_n^k)) = g_i(f_i(a_n^i))$ and the following commutes
\[
\xymatrix{
N_k^* \ar[r]^{g_k} & N^{**} \\
N_{i_0} \ar[u] \ar[r] & N_i^* \ar[u]_{g_i}
}
\]
We put these diagrams together and get the following:
\[
\xymatrix{
 & N^{**} & \\
 N_k^* \ar[ur]_{g_k} & & N_i^* \ar[ul]_{g_i} \\
  & N_{i_0} \ar[ur] \ar[ul] & \\
  N_k \ar[uu]_{f_k} & & N_i \ar[uu]_{f_i} \\
  & M \ar[ur] \ar[ul] \ar[uu] &
}
\]

Thus, we have amalgamated $N^*_i$ and $N^*_k$ over $M$.  Furthermore, for each $j < n+1$, we have $g_k(f_k(a_j^k)) = g_i(f_i(a_j^i))$.  This witnesses $gtp(\seq{a_j^i \mid j < n + 1}/M, N_i) = gtp(\seq{a_j^j \mid j < n + 1}/M, N_j)$, which is a contradiction.\\
Thus, $|gS^{n+1}(M)| \leq \mu = \mathfrak{tb}^\lambda_1$ for all $M \in K_\lambda$ as desired. \hfill \dag\\

This proof can be seen as a semantic generalization of the proof that stability for 1-types implies stability.  Now we wish to prove this upper bound for types of any length $\leq \lambda$.  

The proof works by induction to construct a tree of objects that is indexed by $(\mathfrak{tb}^1_\lambda)$--called $\mu$ in the proof--that codes all $\kappa$ length types as its branches.  Successor stages of the construction are similar to the above proof, but with added bookkeeping.  At limit stages, we wish to continue the construction in a continuous way.  However, we will have a family of embeddings rather than an increasing $\prec_K$-chain.  This is fine since the following closure under direct limits follows from the AEC axioms.

\begin{fact} \label{dirlim}
If we have $\seq{M_i \in K : i < \kappa}$ and, for $i < j < \kappa$, a coherent set of embeddings $f_{i, j}: M_j \to M_i$---that is, one so, for $i < j < k < \kappa$, $f_{i, k} = f_{j, k} \circ f_{i, j}$---then there is an $L(K)$ structure $M = \varinjlim_{i < j < \kappa} (M_i, f_{i, j})$ and embeddings $f_{i, \infty} : M_i \to M$ so that, for all $i < j < \kappa$, $f_{i, \infty} = f_{j, \infty} \circ f_{i, j}$ and, for each $x \in M$, there is some $i < \kappa$ and $m \in M_i$ so $f_{i, \infty}(m) = x$.  Furthermore, the model $M \in K$ and each $f_{i, \infty}$ is a $K$-embedding.
\end{fact}

A proof of this fact can be found in \cite{ramibook}.  This first appeared for AECs in VanDieren's thesis \cite{monicathesis} based on work of Cohn in 1965 on the direct limits of algebras.

We now prove the main theorem.

\begin{theorem} \label{finalone}
If $K$ is an AEC with $\lambda$-AP and $\kappa \leq \lambda$, then $\mathfrak{tb}^\kappa_\lambda \leq (\mathfrak{tb}_\lambda^1)^\kappa$.
\end{theorem}

{\bf Proof:} Set $\mu = \mathfrak{tb}_\lambda^1$.  Let $M \in K_\lambda$ and enumerate $gS^\kappa(M)$ as $\seq{p_i \in gS^\kappa(M) : i < \chi}$, where $\chi = |gS^\kappa(M)|$.  We will show that $\chi \leq \mu^\kappa$, which gives the result.  For each $i < \chi$, find $N_0^i \in K_\lambda$ such that $M \prec N_0^i$ and there is $\seq{a_i^\alpha \in |N_0^i| : \alpha < \kappa} \models p_i$.

The formal construction is laid out below, but we give the idea first.  Our construction will essentially create three objects: a tree of models $\seq{M_\eta : \eta \in {}^{<\kappa} \mu}$; for each $i < \chi$, a function $\eta_i : \kappa \to \mu$; and, for each $i < \chi$, a coherent, continuous system $\{ N^i_\alpha, \bigf^i_{\beta, \alpha} : \beta < \alpha < \kappa \}$.  The tree of models will be domains of types such that the relation of $M_\eta$ to $M_{\eta ^\frown j}$ is like that of $M$ to $N_{i_0}$ in Theorem \ref{firstone}.  We would like the value of the function $\eta_i$ at some $\alpha < \kappa$ to determine the type of $a_i^\alpha$ over $M_{\eta_i \rest \alpha}$.  This can't work because $a_i^\alpha$ isn't in a model also containing $M_\nu$; instead we use its image $\bigf^i_{0, \alpha+1}(a_i^\alpha)$ under the coherent system.  At successor stages of our construction, we will put together elements of equal type over a fixed witness ($i_\eta$ here standing in for $i_0$ in Theorem \ref{firstone}).  At limit stages, we take direct limits.

Once we finish our construction, we show that the map $i \in \chi \mapsto \eta_i \in {}^\kappa \mu$ is injective.  This is done by putting the type realizing sequence together along the chain $\seq{M_{\eta_i \rest \alpha} : \alpha < \kappa}$ to show that $\eta_i$ characterizes $p_i$.

More formally, we construct the following:

\begin{enumerate}
	\item A continuous tree of models $\seq{M_\eta \in K_\lambda : \eta \in {}^{<\kappa}\mu}$ with an enumeration of the types over each model $gS^1(M_\eta) = \{ p_j^\eta : j < |gS^1(M_\eta)| \}$.
	
	\item For each $i < \chi$, a function $\eta_i \in {}^\kappa \mu$.
	
	\item For each $\eta \in {}^{<\kappa} \mu$, an ordinal $i_\eta < \chi$.
	
	\item For each $i < \chi$, a coherent, continuous system $\{N_\alpha^i, \bigf^i_{\beta, \alpha}: N_\beta^i \to_{M_{\eta_i \rest \beta}} N^i_\alpha : \beta < \alpha < \kappa\}$; that is, one such that $\gamma < \beta < \alpha < \kappa$ implies $\bigf^i_{\gamma, \alpha} = \bigf^i_{\beta, \alpha} \circ \bigf^i_{\gamma, \beta}$ and so $\delta < \kappa$ limit implies $(N^i_\delta, \bigf^i_{\alpha, \delta})_{\alpha < \delta} = \varinjlim_{\gamma < \beta < \delta} (N^i_\alpha, \bigf^i_{\gamma, \beta})$.
\end{enumerate}

Our construction will have the following properties for all $\eta \in {}^\beta \mu$ when $\beta < \kappa$.

\begin{enumerate}
	\item[(A)] $i_\eta = \min \{ i < \chi : \eta < \eta_i \}$ if that set is nonempty.

	\item[(B)] $M_{\eta^\frown \seq{j}} := N_\beta^{i_{\eta^\frown \seq{j}}}$ and $M_{\eta_i \rest \beta} \prec N^i_\beta$.
	
	\item[(C)] If $\eta^\frown\seq{j} < \eta_i$, then $p^\eta_j = gtp(\bigf^i_{0, \beta}(a_i^\beta)/M_\eta, N^i_\beta)$.  In particular, this is witnessed by the following diagram
\[
 \xymatrix{
  N_\beta^{i_\eta}  \ar[r]  & N_{\beta+1}^i\\
M_\nu \ar[u] \ar[r] & N_\beta^i \ar[u]_{\bigf_{\beta, \beta+1}^i}
 }
\]	with $\bigf^i_{0, \beta+1}(a_i^\beta) = \bigf^{i_{\eta ^\frown \seq{j}}}_{0, \beta}(a^\beta_{i_{\eta^\frown \seq{j}}})$
\end{enumerate}

{\bf Construction:} At stage $\alpha < \kappa$ of the construction, we will construct $\seq{M_\eta : \eta \in {}^\alpha \mu}$, $\eta_i \rest \alpha$, and $\{ N^i_\alpha, \bigf^i_{\beta, \alpha} : \beta < \alpha\}$ for all $i < \chi$.

$\underline{\alpha = \emptyset}$: We set $M_\emptyset = M$ and note that $N_0^i$ is already defined.  Then $\bigf^i_{0,0}$ is the identity.

$\underline{\alpha\te{ is limit}}$: For each $\eta \in {}^\alpha \mu$, set $M_\eta = \cup_{\beta < \alpha} M_{\eta \rest \beta}$ and $(N^i_\alpha, \bigf^i_{\beta, \alpha})_{\alpha < \delta} = \varinjlim_{\gamma < \beta < \delta} (N^i_\alpha, \bigf^i_{\gamma, \beta})$ as required.  The values of $\eta_i \rest \alpha$ are already determined by the earlier phases of the construction.

$\underline{\alpha = \beta + 1}$: We have constructed our system for each $\nu \in {}^\beta \mu$.  This means that there are enumerations $\{p_k^\nu : k < |gS^1(M_\nu)|\}$ of the 1-types with domain $M_\nu$.  Then, if $i < \chi$ such that $\nu = \eta_i\rest \beta$, we set
$$\eta_i(\beta) = k \te{, where $k < \mu$ is unique such that } gtp(\bigf^i_{0, \beta}(a_i^\beta)/M_\nu, N_\beta^i) = p^\nu_k$$
Then, for each $\eta \in {}^{\alpha} \mu$ set $i_\eta = \min \{ i  < \chi: \eta_i \rest \alpha = \eta \}$ if this set is nonempty; pick it arbitrarily otherwise.  Then, for all $i < \chi$, we have that 
$$gtp(\bigf^i_{0, \beta}(a_i^\beta)/M_\nu, N_\beta^i) = gtp(\bigf^{i_{\eta_i\rest \alpha}}_{0, \beta}(a^\beta_{i_{\eta_i\rest \alpha}})/M_\nu, N_\beta^{{i_{\eta_i\rest \alpha}}})$$
This Galois type equality means that there is a model $N_{\beta+1}^i \succ N_\beta^{i_{\eta_i\rest \alpha}}$ and a function $\bigf^i_{\beta, \beta+1}:N_\beta^i \to_{M_{\nu}} N_{\beta+1}^i$ such that 
$$\bigf^i_{\beta, \beta+1}(\bigf^i_{0, \beta}(a_i^\beta)) = \bigf^{i_{\eta_i\rest \alpha}}_{0, \beta}(a^\beta_{i_{\eta_i\rest \alpha}})$$
Set $M_\eta = N^{\eta_i \rest \alpha}_\beta$ (note that this doesn't depend on the choice of $i$) and, for $\gamma \leq \beta$, set $\bigf^i_{\gamma, \beta+1} = \bigf^i_{\beta, \beta+1} \circ \bigf^i_{\gamma, \beta}$.  This completes the construction.

{\bf This is enough:}  As indicated above, we will show that the map from $i$ to $\eta_i$ is injective.  We do this by showing that $\eta_i = \eta_j$ implies $p_i = p_j$ and, recalling that the enumeration of the $p_i$ were distinct, we must have $i = j$.

Let $i, j < \chi$ such that $\eta := \eta_i = \eta_j$.    We want to show $p_i = p_j$.  We have the following commuting diagram of models for each $\beta < \alpha < \kappa$:
\[
 \xymatrix{
N_0^j \ar[d]_{\bigf_{0,\beta}^j} & M_{\eta \rest 0} \ar[d] \ar[l] \ar[r] & N_0^i \ar[d]_{\bigf_{0, \beta}^i} \\
N_\beta^j \ar[d]_{\bigf_{\beta, \alpha}^j} & M_{\eta \rest \beta} \ar[d] \ar[l] \ar[r] & N_\beta^i \ar[d]_{\bigf_{\beta, \alpha}^i} \\
N_{\alpha}^j  & M_{\eta \rest \alpha} \ar[l] \ar[r] & N_{\alpha}^i  \\
 }
\]
with the property that, for each $\alpha < \kappa$, we know
\begin{eqnarray*}
\bigf_{0,\alpha+1}^i(a_i^\alpha) &=& \bigf_{0, \alpha}^{i_{\eta \rest \alpha+1}} (a_{i_{\eta \rest \alpha+1}}^\alpha)  \\
&=& \bigf_{0,\alpha+1}^j(a_j^\alpha)
\end{eqnarray*}
Note that this element is in $M_{\eta \rest \alpha + 1}$.  Now set $\bigM = \cup_{\alpha< \kappa} M_{\eta \rest \alpha}$.  

Let $k$ stand in for either $i$ or $j$.  Set $(\bigN^k, \bigf^k_{\alpha, \infty})_{\alpha < \kappa} = \varinjlim_{\gamma < \beta < \kappa}(N_\beta^k, \bigf^k_{\gamma, \beta})$.  This gives us the following diagram.
\[
\xymatrix{
N_0^i \ar[d]_{\bigf_{0,\infty}^i} & M \ar[l] \ar[d] \ar[r] & N_0^j \ar[d]_{\bigf_{0, \infty}^j} \\
\bigN^i & \bigM \ar[r] \ar[l] & \bigN^j
}
\]
Then we can amalgamate $\bigN^j$ and $\bigN^i$ over $\bigM$ with
\[
\xymatrix{
\bigN^j \ar[r]^g & N^* \\
\bigM \ar[r] \ar[u] & \bigN^i \ar[u]^f
}
\]
Then, for all $\alpha < \kappa$ and $k = i, j$, $\bigf_{0, \infty}^k (a_k^\alpha) = \bigf_{\alpha+1,\infty}^k( \bigf_{0, \alpha+1}^k (a_k^\alpha))$.  We know that $\bigf_{0, \alpha+1}^k (a_k^\alpha) \in |M_{\eta \rest \alpha+1}|$, so it is fixed by $f_\beta^k$ for $\beta > \alpha +1$.  This means it is also fixed by $\bigf_{\alpha+1,\infty}^k$.  Then 
$$\bigf_{0,\infty}^k (a_k^\alpha) = \bigf_{\alpha+1,\infty}^k( \bigf_{0, \alpha+1}^k (a_k^\alpha) ) =  \bigf_{0, \alpha+1}^k (a_k^\alpha) = \bigf_{0, \alpha}^{i_{\eta \rest \alpha+1}} (a_{i_{\eta \rest \alpha+1}}^\alpha) $$
Since this last term is independent of whether $k$ is $i$ or $j$, we have $\bigf_{0,\infty}^i(a_i^\alpha) = \bigf_{0,\infty}^j(a_j^\alpha) \in \bigM$ for all $\alpha < \kappa$.  Since our amalgamating diagram commutes over $\bigM$, $f(\bigf_0^i(a_i^\alpha)) = g(\bigf_0^j(a_j^\alpha))$.\\
Combining the above, we have
\[
\xymatrix{
N_0^j \ar[r]^{g \circ \bigf_{0,\infty}^j} & N^* \\
M \ar[u] \ar[r] & N_0^i \ar[u]_{f \circ \bigf_{0,\infty}^i}
}
\]
with $f \circ \bigf_{0,\infty}^i(\seq{a_i^\alpha \mid \alpha < \kappa}) = g \circ \bigf_{0,\infty}^j (\seq{a_j^\alpha \mid \alpha < \kappa})$.\\
Thus,
$$p_i = gtp(\seq{a_i^\alpha \mid \alpha < \kappa} / M, N_0^i) = gtp(\seq{a_j^\alpha \mid \alpha < \kappa} / M, N_0^j) = p_j$$
Since each $p_k$ was distinct, this implies that $i=j$.  The map $i \mapsto \eta_i$ is injective and $\chi \leq \mu^\kappa$ as desired. \hfill \dag\\

As mentioned in the Introduction, the above result gives us the proof of Theorem \ref{introtheorem}:

{\bf Proof of Theorem \ref{introtheorem}:} As discussed in the last section, $(\te{Mod }T, \prec_{L(T)})$ is an AEC with amalgamation over sets.  Given a set $A$, passing to a model containing $A$ can only increase the number of types.  Thus, even in this case, it is enough to only consider models when computing $\mathfrak{tb}$. Thus, 
$$\sup_{A \subset M\models T, \|A\| = \lambda} |S^\mu(A)| = \mathfrak{tb}_\lambda^\mu = (\mathfrak{tb}_\lambda^1)^\mu = \Bigg(\sup_{A \subset M\models T, \|A\| = \lambda} |S^1(A)| \Bigg)^\mu$$
as desired.  \hfill \dag\\

After seeing this work, Alexei Kolesnikov pointed out a much simpler proof of Theorem \ref{finalone} for first order theories or, more generally, for AECs that are $< \omega$ type short over $\lambda$-sized domains\footnote{See \cite{tamelc}.3.4 for a definition or ignore this case at no real loss}; in either case, a type of infinite length is determined by its restrictions to finite sets of variables.  Fix a type $p \in S^I(M)$ with $I$ infinite.  The previous comment means that the map
$$p \mapsto \Pi_{\bx \in [I]^{<\omega}} p^\bx$$
from $S^I(M)$ to $\Pi_{\bx \in [I]^{<\omega}} S^\bx(M)$ is injective.  Then
\begin{eqnarray*}
|S^I(M)| &\leq& \Pi_{\bx \in [I]^{<\omega}} |S^\bx(M)| = \Pi_{\bx \in [I]^{<\omega}} |S^1(M)|\\
&=& |S^1(M)|^{|[I]^{<\omega}|} = |S^1(M)|^{|I|}
\end{eqnarray*}
This is in fact a strengthening of Theorem \ref{finalone} as in Theorem \ref{lowerbound}.

We now examine local types in first order theories.  For $\Delta \subset Fml(L(T))$, set
$$\Delta\mathfrak{tb}_\kappa^\lambda = \sup_{M \models T, \|M\| = \lambda} |S^\kappa_\Delta(M)|$$
If $\Delta = \{\phi\}$, we simply write $\phi\mathfrak{tb}^\kappa_\lambda$.  Unfortunately, there is no semantic equivalent of $\Delta$-types, so the methods and proofs above do not transfer.  For a lower bound, we can prove the following in the same way as Theorem \ref{lowerbound}.

\begin{prop} \label{locallowerbound}
If $T$ is a first order theory and $\Delta \subset Fml(L(T))$, then for any $\kappa$ we have that $|S^1_\Delta(A)| = \mu$ implies that $|S^\kappa_\Delta(A)| \geq \mu^\kappa$.
\end{prop}

If $\Delta$ is closed under existential quantification, the syntactic proofs of Theorem \ref{firstone} (see, for instance, \cite{shelahfobook}.2) can be used to get an upper bound for $\Delta\mathfrak{tb}^n_\lambda$ when $n$ is finite.

\begin{prop}\label{3point7}
If for all $\phi(\bx, x, \by) \in \Delta$, we have $\exists z \phi(\bx, z, \by) \in \Delta$, then $\Delta\mathfrak{tb}_\lambda^n \leq \Delta\mathfrak{tb}_\lambda^1$ for $n<\omega$.
\end{prop}

With this result for finite lengths, we can apply the syntactic argument above to conclude the following.

\begin{prop}
If $T$ is a first order theory and $\Delta \subset Fml(L(T))$, then for $\kappa \leq \lambda$, 
$$\Delta\mathfrak{tb}^\kappa_\lambda \leq (\sup_{n < \omega} \Delta\mathfrak{tb}^n_\lambda)^\kappa$$
In particular, if $\Delta$ is closed under existentials as in Proposition \ref{3point7}, then $\Delta\mathfrak{tb}^\kappa_\lambda \leq (\Delta\mathfrak{tb}^1_\lambda)^\kappa$.
\end{prop}

We now turn to the values of $\phi\mathfrak{tb}$ for particular $\phi$.  Recall that Theorem \cite{shelahfobook}.II.2.2 says that $T$ is stable iff $T$ is $\lambda$ stable for $\lambda = \lambda^{|T|}$ iff it is $\lambda$ stable for $\phi$ types for all $\phi \in L(T)$.  This means that if $T$ is unstable in $\lambda = \lambda^{|T|}$, then there is some $\phi$ such that $\phi\mathfrak{tb}_\lambda^1 > \lambda$.
Further, suppose that $\sup \{ \psi\mathfrak{tb}_\lambda^1: \psi \in L(T) \} = \lambda^{+n}$ for some $1 \leq n < \omega$.  Then, since $\lambda^{+n}$ is a successor, this supremum is acheived by some formula $\phi_\lambda$.  Then, since $\lambda^{|T|} = \lambda$, we can calculate
\begin{eqnarray*}
\phi_\lambda \mathfrak{tb}^1_\lambda &=& \sup_{\psi \in L(T)} \{ \psi  \mathfrak{tb}_\lambda^1 \} \leq \mathfrak{tb}_\lambda^1 \leq \Pi_{\psi \in L(T)} (\psi\mathfrak{tb}_\lambda^1) \leq (\phi_\lambda\mathfrak{tb}_\lambda^1)^{|T|} = \\
&=& (\lambda^{+n})^{|T|} = \lambda^{|T|} \cdot \lambda^{+n} = \lambda^{+n} = \phi_\lambda \mathfrak{tb}_\lambda^1
\end{eqnarray*}
So $\phi_\lambda\mathfrak{tb}_\lambda^1 = \mathfrak{tb}_\lambda^1$.  Thus, for all $\kappa \leq \lambda$, we can use Theorems \ref{locallowerbound} and \ref{finalone} to calculate
$$(\phi_\lambda\mathfrak{tb}_\lambda^1)^\kappa \leq \phi_\lambda\mathfrak{tb}_\lambda^\kappa \leq \mathfrak{tb}_\lambda^\kappa = (\mathfrak{tb}_\lambda^1)^\kappa = (\phi_\lambda\mathfrak{tb}_\lambda^1)^\kappa$$
This gives us the following result:

\begin{theorem}
Given a first order theory $T$, if $\lambda$ is a cardinal such that $\lambda^{|T|} = \lambda$ and $\sup \{ |S^1_\psi(A)| : \psi \in L(T), |A|\leq \lambda\} < \lambda^{+\omega}$, then there is some $\phi_\lambda \in L(T)$ such that, for all $\kappa \leq \lambda$, $\mathfrak{tb}_\lambda^\kappa = \phi_\lambda\mathfrak{tb}_\lambda^\kappa$.
\end{theorem}

\label{modelextension}

Returning to general AECs, in \cite{shelahaecbook}.V.D.\S3, Shelah considers long types of tuples enumerating a model extending the domain.  In this case, any realization of the type is another model extending the domain that is isomorphic to the original tuple over the domain.  Thus, an upper bound on types of a certain length $\kappa$ also bounds the number of isomorphism classes extending the domain by $\kappa$ many elements.  More formally, we get
\begin{prop}
Given $M \in K_\lambda$,
$$|\{ N/\cong_M : N \in K, M \precneqq N, |N - M| = \kappa \}| \leq \mathfrak{tb}_\lambda^\kappa$$
\end{prop}
If we have an AEC with amalgamation where any extension can be broken into smaller extensions, this could lead to a useful analysis.  Unfortunately, this provides us with no new information when $\kappa = \lambda$ since $2^\lambda = \mathfrak{tb}_\lambda^\lambda$ is already the well-known upper bound for $\lambda$-sized extensions of $M$ and there are even first order theories where $M \precneqq N$ implies $|N - M| \geq \| M \|$.  Algebraically closed fields of characteristic 0 are such an example.

\section{Strongly Separative Types} \label{strongseptypes} 

One might hope that similar bounds could be developed for non-algebraic types.  This would probably give us a finer picture of what is going on because a model $M$ necessarily has at least $\|M\|$ many algebraic types over $M$, so in the stable case, the number of non-algebraic types could, {\it a priori}, be anywhere between 0 and $\|M\|$; the case $gS^1_{na}(M) = \emptyset$ only occurs in the uninteresting case that $M$ has no extensions.

However, as the following result shows, no such result is possible even in basic, well-understood first order cases:

\begin{prop} \label{badnaex}
\begin{enumerate}
	\item Let $T_1$ be the empty theory and $M \vDash T_1$.  Then $|S^1_{na}(M)| = 1$ and $|S^n_{na}(M)| = B_n$ for all $n < \omega$, where $B_n$ is the $n$th Bell number.  In particular, this is finite.
	
	\item Let $T_2 = ACF_0$ and $M \vDash T$.  Then $|S^1_{na}(M)| = 1$ but $|S^2_{na}(M)| = \|M\|$.
	
\end{enumerate}
\end{prop}

Note that these examples represent the minimal and maximal, respectively, number of long, nonalgebraic types given that there is only one non-algebraic type.

{\bf Proof:} \begin{enumerate}

	\item Let $tp(a/M, N_1), tp(b/M, N_2) \in S^1_{na}(M)$ and, WLOG, assume $\|N_1\| \leq \|N_2\|$.  Then let $f$ fix $M$, send $a$ to $b$, and injectively map $N_1 - M - \{a\}$ to $N_2 - M - \{b\}$ arbitrarily.  This witenesses $tp(a/M, N_1) = tp(b/M, N_2)$.\\
	Given the type of $\seq{a_n : n < k}$, the only restriction on finding a function to witness type equality is given by which elements of the sequence are repeated; for instance, if $a \neq b$, then $tp(a, b/M, N) \neq tp(a, a/M, N)$.  Thus, each type can be represented by those elements of the sequence which are repeated.  To count this, we need to know the number of partitions of $n$.  This is given by Bell's numbers, defined by $B_1 = 1$ and $B_{n+1} = \sum_{k=0}^n {n \choose k} B_k$. See \cite{generatingfunctionology}.1.6.13 for a reference.  Then, the number of $n$-types is just $B_n$.
	
	\item This is an easy consequence of Steinitz's Theorem that there is only one non-algebraic 1-type, that of an element transcedental over the domain.\\
	Looking at 2-types, let $M \in K$ and $e \in N \succ M$ be transcendental (non-algebraic) over $M$.  For each polynomial $f \in M[x]$, set $p_f = tp(e, f(e)/M, N)$.  Then, for $f \neq g$, we have that $p_f \neq p_g \in S_{na}^2(M)$. Also, $tp(e, \pi/M, N')$ is distinct, where $\pi$ is transcendental over $M(e)$.  This gives at least $\|M\|$ many 2-types.\\
	We know there are at most $\|M\|$ many because the theory is stable.  Therefore, the results of last section tells us that there are exactly $\|M\|^2 = \|M\|$ many 2-types, so there are at most $\|M\|$ many non-algebraic 2-types.\hfill \dag\\
	
\end{enumerate}

This shows that a result like Theorem \ref{finalone} is impossible for non-algebraic types.  As is evident in the proof above, especially part two, the variance in the number of types comes from the fact that, while the realizations of the non-algebraic type are not algebraic over the model, they might be algebraic over each other.  This means that even 2-types, like $tp(e, 2e/\mathbb{A}, \mathbb{C})$, that are not realized in the base model can't be separated: any algebraically closed field realizing the type of $e$ must also realize the type of $2e$.

In order to get a bound on the number of these types, we want to be able to {\it separate} the different elements of the tuples that realize the long types.  This motivates our definition and naming of separative types below.  We also introduce a slightly stronger notion, strongly separative types, that allow us to not only separate realizations of the type, but also gives us the ability to extend types, as made evident in Proposition \ref{strsepclosure}.  Luckily, in the first order case and others, these two notions coincide; see Proposition \ref{dapsep}.

\begin{defin} \label{sepdef}
\begin{enumerate}

	\item[]
	
	\item We say that a triple $( \seq{a_i : i < \alpha}, M, N) \in K^{3, \alpha}_\lambda$ is \emph{separative} iff there are increasing sequences of intermediate models $\seq{N_i \in K : i < \alpha}$ such that, for all $i < \alpha$, $M \prec N_i \prec N$ and $a_i \in N_{i+1} - N_i$.  The sequence $\seq{N_i : i < \alpha}$ is said to witness the triple's separativity.

	\item For $M \in K$, set $gS^\alpha_{sep}(M) = \{ gtp(\seq{a_i : i < \alpha,}/ M, N) : (\seq{a_i : i < \alpha,}, M, N) \in K_\lambda^{3, \alpha}$ is separative$ \}$.

	\item We say that a triple $( \seq{a_i : i < \alpha}, M, N) \in K^{3, \alpha}_\lambda$ is \emph{strongly separative} iff there is a sequence witnessing its separativity $\seq{N_i : i < \alpha}$ that further has the property that, for any $i < \alpha$ and $N_1^+ \succ N_i$ of size $\lambda$, there is some $N_2^+ \succ N_1^+$ and $g: N_{i+1} \to_{N_\beta} N_2^+$ such that $g(a_i) \notin N_1^+$.
	
	\item For $M \in K_\lambda$, set $gS^\alpha_{strsep}(M) = \{ gtp(\seq{a_i : i < \alpha,}/ M, N) : (\seq{a_i : i < \alpha,}, M, N) \in K_\lambda^{3, \alpha}$ is strongly separative$\}$.	

\end{enumerate}
\end{defin}

The condition ``$a_i \in N_{i+1} - N_i$'' in (1) could be equivalently stated as either of the following:
\begin{itemize}
	\item For all $j < \alpha$, $a_j \in N_i$ iff $i < j$.
	\item $gtp(a_i/N_i, N_{i+1})$ is nonalgebraic.
\end{itemize}

Note that the examples in Proposition \ref{badnaex} only have one separative or strongly separative type of any length: for the empty theory, this is any sequence of distinct elements and, for $ACF_0$, this is any sequence of mutually transcendental elements.  Theorem \ref{strsepone} below shows this generally by proving the upper bound from the last section (Theorem \ref{finalone}) holds for strongly separative types.  Before this proof, a few comments about these definitions are in order.

First, the key part of the definition is about triples, but we will prove things about types.  This is not an issue because any triple realizing a (strongly) separative type can be made into a (strongly) separative type by extending the ambient model.

\begin{prop}
\begin{enumerate}

	\item[]
	
	\item If $gtp(\seq{a_\beta: \beta < \alpha}/M, N) \in gS^\alpha_{sep}(M)$, then there is some $N^+ \succ N$ such that $(\seq{a_\beta: \beta < \alpha}, M, N^+)$ is separative.
	
	\item The same is true for strongly separative types.

\end{enumerate}
\end{prop}

{\bf Proof:}  We will prove the first assertion and the second one follows similarly.  By the definition of $gS^\alpha_{sep}$, there is some separative $(\seq{b_\beta : \beta < \alpha}, M, N_1) \in K^{3, \alpha}_\lambda$ such that $gtp(\seq{a_\beta: \beta < \alpha}/M, N) = gtp(\seq{b_\beta: \beta < \alpha}/M, N)$.  Thus, there exists some $N^+ \succ N$ and $f: N_1 \to_M N^+$ such that $f(b_\beta) = a_\beta$ for all $\beta < \alpha$.  Let $\seq{N_\beta:\beta < \alpha}$ be a witness sequence to $(\seq{b_\beta : \beta < \alpha}, M, N_1)$'s separativity.  Then $\seq{f(N_\beta) \prec N^+:\beta < \alpha}$ is a witness sequence for $(\seq{a_\beta : \beta < \alpha}, M, N^+)$.\hfill \dag\\

Second, although we continue to use the semantic notion of types (Galois types) for full generality, these notions are new in the context of first order theories.  In this context, the elements of the witnessing sequence $\seq{N_i : i < \alpha}$ are still required to be models, even though types are meaningful over sets.  An attempt to characterize these definitions in a purely syntactical nature (i.e. by only mentioning formulas) was unsuccessful, but we do know (see Proposition \ref{dapsep} below) that all separative types over models are strongly separative for complete first order theories.

Third, we can easily characterize these properties for 1-types.

\begin{prop} \label{1typechar}
Let $K$ be an AEC and $p \in gS^1(M)$.
\begin{itemize}

	\item $p$ is separative iff $p$ is nonalgebraic.
	
	\item $p$ is strongly separative iff, for any $N \succ M$ with $\|N\| = \|M\|$, there is an extension of $p$ to a non-algebraic type over $N$.  Such types are called \emph{big}.

\end{itemize}
\end{prop}

Finally, strongly separative types and separative types are the same in the presence of the disjoint amalgamation property.

\begin{prop} \label{dapsep}
Let $\alpha$ be an ordinal and $M \in K$.  If $K$ satisfies the disjoint amalgamation property when all models involved have sizes between $\|M\|$ and $|\alpha| + \|M\|$, inclusive, then $gS^\alpha_{strsep}(M) = gS^\alpha_{sep}(M)$.
\end{prop}

{\bf Proof:}  By definition, $gS^\alpha_{strsep}(M) \subset gS^\alpha_{sep}(M)$, so we wish to show the other containment.  Let $gtp(\seq{a_\beta : \beta < \alpha}/M, N) \in gS^\alpha_{sep}(M)$.  Let $\seq{N_\beta: \beta<\alpha}$ be a witnessing sequence and let $N_1^+ \succ N_{\beta_0}$ of size $\|N_{\beta_0}\|$ for some $\beta_0 < \alpha$.  By renaming elements, we can find some copy of $N_1^+$ that is disjoint from $N_{\beta_0+1}$ except for $N_{\beta_0}$.  So there are $\bigN$ and $f:N_1^+ \cong_{N_{\beta_0}} \bigN$ such that $\bigN \cap N_{\beta_0+1} = N_{\beta_0}$.  Then, we can use disjoint amalgamation on $\bigN$ and $N_{\beta_0+1}$ over $N_{\beta_0}$ to get $N^*$ and $g:N_{\beta_0+1} \to N^*$ so
\[
\xymatrix{
N_1^+ \ar[r]^f & \bigN \ar[r] & N^*\\
& N_{\beta_0} \ar[ul] \ar[u] \ar[r] & N_{\beta_0+1} \ar[u]_g
}
\]
commutes and $\bigN \cap g(N_{\beta_0+1}) = N_{\beta_0}$.  Thus, since $a_{\beta_0}$ is in $N_{\beta_0+1}$ and not in $N_{\beta_0}$, we have that $g(a_{\beta_0})$ is in $g(N_{\beta_0+1})$ and not in $\bigN$.  Let $\widehat{f}$ be an $L(K)$-isomorphism that extends $f$ and has $N^*$ in its range.  Then we have
$$\widehat{f}^{-1}(g(a_{\beta_0})) \not \in \widehat{f}^{-1}(\bigN) = f^{-1}(\bigN) = N_1^+$$
Then we can collapse the above diagram to 
\[
\xymatrix{
N_1^+ \ar[r] & \widehat{f}^{-1}(N^*) \\
N_{\beta_0} \ar[u] \ar[r] & N_{\beta_0+1} \ar[u]_{\widehat{f}^{-1} \circ g}
}
\]
This diagram commutes and witnesses the property for strong separativity with $N^+_2 = \widehat{f}^{-1}(N^*)$.  \hfill \dag\\

It is an exercise in the use of compactness that every complete first order theory satisfies disjoint amalgamation over models, see Hodges \cite{modeltheoryhodges}.6.4.3 for a reference.  For a general AEC, this is not the case.  Baldwin, Kolesnikov, and Shelah \cite{bks} have constructed examples of AECs without disjoint amalgamation.  On the other hand, Shelah \cite{sh576} has shown that disjoint amalgamation follows from certain amounts of structure (see, in particular, 2.17 and 5.11 there).  Additionally, Grossberg, VanDieren, and Villaveces \cite{gvv} point out that many AECs with a well developed independence notion, such as homogeneous model theory or finitary AECs, also satisfy disjoint amalgamation.

In order to prove the main theorem of this section, Theorem \ref{strsepone}, we will need to make use of certain closure properties of strongly separative types.  These also hold for separative types as well.

\begin{prop}[Closure of $gS_{strsep}$] \label{strsepclosure}
\begin{enumerate}

	\item[]

	\item If $p \in gS^\alpha_{strsep}(M)$ and $I \subset \alpha$, then $p^I \in gS^{otp(I)}_{strsep}(M)$.
	
	\item If $p \in gS^\alpha_{strsep}(M)$ and $M_0 \prec M$, then $p \rest M_0 \in gS^\alpha_{strsep}(M_0)$.
	
\end{enumerate}
\end{prop}

We now prove the main theorem.

\begin{defin}
 The \emph{strongly separative type bound for $\lambda$ sized domains and $\kappa$ lengths} is denoted $strsep\mathfrak{tb}^\kappa_\lambda = \sup_{M \in K_\lambda} |gS^\kappa_{strsep}(M)|$.
\end{defin}

\begin{theorem} \label{strsepone}
If $strsep\mathfrak{tb}_\lambda^1 = \mu$, then $strsep\mathfrak{tb}_\lambda^\kappa \leq \mu^\kappa$ for all (possibly finite) $\kappa \leq \lambda^+$.
\end{theorem}

{\bf Proof:}  The proof is very similar to that of Theorem \ref{finalone}, so we only highlight the differences.

As before, let $M \in K_\lambda$, enumerate $gS^\kappa_{strsep}(M) = \seq{p_i : i < \chi}$, and find $N_0^i \succ M$ of size $\lambda + \kappa$ and $a_i^\alpha \in |N_0^i|$ for $i < \chi$ and $\alpha < \kappa$ such that $\seq{a_i^\alpha : \alpha < \kappa} \vDash p_i$.

Then, we use strong separativity to find a witnessing sequence.  That is, for each $i < \chi$, we have increasing and continuous $\seq{{}^\alpha N_0^i \in K_\lambda : \alpha < \kappa}$ so, for each $\alpha < \kappa$, $M \prec {}^\alpha N^i_0 \prec N^i_0$ and $a^\alpha_i \in {}^{\alpha + 1}N_0^i - {}^\alpha N_0^i$.

As before, we will construct $\seq{M_\eta \in K_\lambda : \eta \in {}^{< \kappa} \mu}$, $\seq{p_j^\eta \in gS^1_{strsep}(M_\eta) : j < |gS^1_{strsep}(M_\eta)|}$, $\seq{i_\eta \in \chi : \eta \in {}^{<\kappa}\mu}$, and $\seq{\eta_i \in {}^\kappa \mu : i < \chi}$ as in $(1)-(3)$ of the proof of Theorem \ref{finalone} and
\begin{enumerate}

	\item[(4*)] For $i < \chi$, a coherent, continuous $\{{}^\alpha N_\alpha^i, \bigf^i_{\beta, \alpha}: {}^\beta N^i_\beta \to_{M_{\eta_i \rest \beta}} {}^\alpha N_\alpha^i \mid \beta < \alpha < \kappa \}$, models $\seq{ {}^{\alpha+1} N_\alpha^i : \alpha < \kappa}$, and, for each $\beta < \alpha < \kappa$, functions
	\begin{itemize}
		\item $h_\alpha^i: {}^\alpha N_0^i \to {}^\alpha N_\alpha^i$;
		\item $g^i_{\beta+1}: {}^{\beta+1} N_0^i \to {}^{\beta+1} N_\beta^i$; and 
		\item $f^i_{\beta+1}: {}^{\beta+1} N_\beta^i \to_{M_{\eta_i \rest \beta}} {}^\beta N_\beta^i$.
	\end{itemize}

\end{enumerate}	
These will satisfy (A), (B), and (C) from Theorem \ref{finalone} and
\begin{enumerate}	
	
	\item[(D)] If $\alpha = \beta+1$, then $h^i_\alpha = f^i_\alpha \circ g^i_\alpha$ and, if $\alpha$ is limit, then ${}^\alpha N_\alpha^i$ is the direct limit and, for each $\delta < \alpha$, the following commutes
			\[
		 \xymatrix{\ar @{} [dr] {}^\delta N_\delta^i  \ar[r]^{\bigf^i_{\delta, \alpha}}  & {}^\alpha N_\alpha^i \\
		{}^\delta N_0^i \ar[u]^{h_\delta^i} \ar[r] & {}^\alpha N_0^i \ar[u]_{h_\alpha^i}
		 }
		\]

	\item[(E)] If $\alpha < \kappa$, then
	\begin{itemize}
		\item $h^i_\alpha \rest {}^\alpha N^i_0 = g^i_{\alpha+1} \rest {}^\alpha N^i_0$
		\item $g^i_{\alpha+1}(a_i^\alpha) \not\in {}^\alpha N_\alpha^i$
		\item $h^i_{\alpha+1}(a_i^\alpha) = g^{i_{\eta_i\rest \alpha}}_{\alpha+1}(a^\alpha_{i_{\eta_i \rest \alpha}})$
		\item $\bigf^i_{\alpha, \alpha+1} = f^i_{\alpha+1}$
	\end{itemize}

\end{enumerate}

{\bf Construction:}

The base case and limit case are the same as in \ref{finalone}.  In the limit, we additionally set $h_\alpha^i = \bigcup_{\beta<\alpha} h^i_\beta $.

For $\ell(\eta) = \alpha = \beta + 1$ we will apply our previous construction to the separating models.  Fix some $\nu \in {}^\beta \mu$.  For each $i < \chi$ such that $\eta_i \rest \beta = \nu$, we have $h^i_\beta: {}^\beta N_0^i \to {}^\beta N_\beta^i$.  We know that $gtp(a_i^\beta/{}^\beta N^i_0, {}^{\beta+1}N^i_0)$ is big by Propositions \ref{strsepclosure} and \ref{1typechar}.  Thus we can find a big extension with domain $(h^i_\beta)^{-1}({}^\beta N_\beta^i)$.  Then, applying $h^i_\beta$ to this type, we get some $g^i_{\beta+1} : {}^{\beta+1} N_0^i \to {}^{\beta+1} N_\beta^i$ so
	\[
	 \xymatrix{\ar @{} [dr] 
	 {}^\beta N_\beta^i  \ar[r]  & {}^{\beta+1}N_\beta^i\\
	{}^\beta N_0^i \ar[u]^{h_\beta^i} \ar[r] & {}^{\beta+1}N_0^i \ar[u]_{g_{\beta+1}^i}
	 }
	\]
	commutes and $gtp(g_{\beta+1}^i(a_i^\beta)/M_\nu, {}^{\beta+1} N^i_\beta)$ is big and, therefore, strongly separative.  Note that this extension uses that these types are strongly separative and not just separative.  Then we can extend $\eta_i$ by $\eta_i(\beta) = k$ where $k < \mu$ is the unique index such that $gtp(g_{\beta+1}^i(a_i^\beta)/M_\nu, {}^{\beta+1} N^i_\beta) = p^\nu_k$.\\
	Then set $i_{\nu^\frown \seq{i}} = \min \{ i < \chi : \eta_i \rest \alpha = \nu^\frown \seq{i}\}$.  This means that, for all $i <\chi$, we have 
	$$gtp(g_{\beta+1}^i(a_i^\beta)/M_\nu, {}^{\beta+1}N_\beta^i) = gtp(g_{\beta+1}^{i_{\eta_i \rest \alpha}}(a_{i_{\eta_i \rest \alpha}}^\beta)/M_\nu, {}^{\beta+1}N_\beta^{i_{\eta_i \rest \alpha}})$$
	Thus, we can find ${}^{\beta+1} N_{\beta+1}^i \succ {}^{\beta+1}N_\beta^{i_{\eta_i \rest \alpha}}$ from $K_\lambda$ and $f_{\beta+1}^i:{}^{\beta+1}N_\beta^i \to_{M_\nu} {}^{\beta+1}N_{\beta+1}^i$ such that $f_{\beta+1}^i(g_{\beta+1}^i(a_i^\beta)) = g_{\beta+1}^{i_{\eta_i \rest \alpha}}(a_{i_{\eta_i \rest \alpha}}^\beta)$.  Finally, set $M_{i_{\eta_i \rest \alpha}} = {}^{\beta+1}N_\beta^{i_{\eta_i \rest \alpha}}$ and $h^i_{\beta+1} = f^i_{\beta+1} \circ g^i_{\beta+1}$.

{\bf This is enough:}  For each $i < \chi$ and every $\alpha < \beta < \kappa$, we have
\[
\xymatrix{
M_{\eta\rest 0} \ar[r] \ar[d] & {}^0 N_0^i \ar[r] \ar[d]_{\bigf^i_{0, \beta}} & {}^\beta N_0^i \ar[r] \ar[dl]^{h_\beta^i} & {}^\alpha N_0^i \ar[r] \ar[ddll]^{h_\alpha^i} & N_0^i \\
M_{\eta \rest \beta} \ar[d] \ar[r] & {}^\beta N_\beta^i \ar[d]_{\bigf_{\beta, \alpha}^i} \\
M_{\eta \rest \alpha} \ar[r] & {}^\alpha N_\alpha^i
}
\]
that commutes.  Note that this is almost the same diagram as before, except we have added the separating sequences.  Then we can proceed as before, setting
\begin{enumerate}

	\item $\bigM = \cup_{\alpha < \kappa} M_{\eta \rest \alpha}$;
	
	\item $(\bigN^i , \bigf_{\alpha, \infty}) = \varinjlim_{\beta < \gamma < \kappa} ({}^\beta N_{\beta}^i, \bigf^i_{\alpha, \beta})$;
	
	\item $N_1^i = \cup_{\alpha < \kappa} {}^\alpha N_0^i \prec N_0^i$.
	
	\item $f_i: N_1^i \to \bigN^i$ by $f_i = \cup_{\alpha < \kappa} ( \bigf_{\alpha, \infty} \circ h_\alpha^i)$; and
	
	\item $\eta_i \in {}^\kappa \mu$ such that $i \in I_{\eta \rest \alpha}$ for all $\alpha < \kappa$. 
	
\end{enumerate}

Then, if $\chi > \mu^\kappa$, there are $i \neq j$ such that $\eta_i = \eta_j$.  As before, this would imply $p_i = p_j$, but they are all distinct.  So $\chi \leq \mu^\kappa$ as desired.   \hfill \dag\\

In the previous theorem, we allowed the case $\kappa = \lambda^+$.  Most of the time, this is only the set-theoretic bound $strsep\mathfrak{tb}_\lambda^{\lambda^+} \leq 2^{\lambda^+}$.  However, if we had $strsep\mathfrak{tb}_\lambda^1 = 1$, then we get the surprising result that $strsep\mathfrak{tb}_\lambda^{\lambda^+} \leq 1$.  This will be explored along with further investigation of classifying AECs based on separative types in future work.

\section{Saturation} \label{satsec}

We now turn from the number of infinite types to their realizations.  The saturation version of Theorem \ref{aectheorem} is much simpler to prove.

\begin{prop} \label{basicsat}
If $M \in K_\lambda$ is Galois saturated for 1-types, then $M$ is Galois saturated for $\lambda$-types.
\end{prop}

{\bf Proof:} Let $M_0 \prec M$ of size $< \lambda$ and $p \in gS^\lambda(M_0)$.  By the definition of Galois types, there is some $N \succ M_0$ of size $\lambda$ that realizes $p$.  Find a resolution of $N$ $\seq{N_i \in K_{< \lambda} \mid i < \cf \lambda}$ with $N_0 = M_0$.  Then use Lemma \ref{modhomsat} to get increasing, continuous $f_i:N_i \to M$ that fix $M_0$.  Then $f := \cup_{i < \lambda} f_i : N \to_{M_0} M$.  This implies $f(N) \models f(p) = p$ and since $f(N) \prec M$, $M \models p$. \hfill \dag

We can get a parameterized version with the same proof.

\begin{prop}
If $M \in K_\lambda$ is $\mu$-Galois saturated for 1-types, then $M$ is $\mu$-Galois saturated for $\mu$-types.
\end{prop}

The seeming simplicity of the proof of Proposition \ref{basicsat}, especially compared to earlier uses of direct limits, hides the difficulty and complexity of the proof of Lemma \ref{modhomsat}.  Although the statement is a generalization of a first order fact, its announcement was a surprise and many flawed proofs were proposed before a successful proof was given.  Building on work of Shelah, Grossberg, and Kolesnikov, Baldwin \cite{baldwinbook}.16.5 proves a version of Lemma \ref{modhomsat} which does not require amalgamation.  This gives rise to a version of Proposition \ref{basicsat} in AECs even without amalgamation.

There is also a strong relationship between the value of $\mathfrak{tb}_\lambda^1$ and the existence of $\lambda^+$-saturated extensions of models of size $\lambda$.  The following generalizes first order theorems like \cite{shelahfobook} Theorem VIII.4.7.  

In the following theorems, we make use of a monster model, as in first order model theory, to reduce the complexity of constructions.  Full details can be found in the references given at the start of Section \ref{prelimsec}, but the key facts are
\begin{itemize}
	\item the existence of a monster model $\sea$ follows from the amalgamation property, the joint embedding property, and every model having a proper $\prec_K$-extension; and
	\item for $M \prec \sea$ and $a, b \in |\sea|$, 
	$$gtp(a/M) = gtp(b/M) \iff \exists f \in Aut_M \sea \te{ so } f(a) = b$$
\end{itemize}

The first relationship is clear from counting types.  

\begin{theorem}
Let $K$ be an AEC with amalgamation, joint embedding, and no maximal models. If every $M \in K_\kappa$ has an extension $N \in K_\lambda$ that is $\kappa^+$-saturated, then $\mathfrak{tb}_\kappa^1 \leq \lambda$.
\end{theorem}

{\bf Proof:} Assume that every model in $K_\kappa$ has a $\kappa^+$-saturated extension of size $\lambda$.  Let $M \in K_\kappa$ and $N \in K_\lambda$ be that extension.  Since every type over $M$ is realized in $N$, we have $|gS(M)| \leq \|N\| = \lambda$.  Taking the sup over all $M \in K_\kappa$, we get $\mathfrak{tb}_\kappa^1 \leq \lambda$, as desired. \hfill \dag

Going the other way, we have both a set theoretic hypothesis and model theoretic hypothesis that imply instances of a $\kappa^+$-saturated extension.  The set-theoretic version is well known.

\begin{theorem}
Let $K$ be an AEC with amalgamation, joint embedding, and no maximal models.  If $\lambda^\kappa = \lambda$, then every $M \in K_\kappa$ has an extension $N \in K_\lambda$ that is $\kappa^+$ saturated.
\end{theorem}

Note that the hypothesis implies $\mathfrak{tb}_\lambda^1 \leq \lambda$.  Without this set theoretic hypothesis, reaching our desired conclusion is much harder.  $\lambda^\kappa = \lambda$ means that we can consider all $\kappa$ size submodels of a $\lambda$ sized model without going up in size.  Without this assumption, things become much more difficult and we must rely on model theoretic hypotheses.  The following has a stability-like hypothesis, sometimes called `weak stability;' see \cite{jrsh875}, for instance.

\begin{theorem}
Let $K$ be an AEC with amalgamation, joint embedding, and no maximal models.  If $\mathfrak{tb}_\kappa^1 \leq \kappa^+$, then every $M \in K_\kappa$ has an extension $N \in K_{\kappa^+}$ that is saturated.
\end{theorem}

{\bf Proof:} We proceed by a series of increasingly strong constructions.\\
{\bf Construction 1:} For all $M \in K_\kappa$, there is $M^* \in K_{\kappa^+}$ such that all of $S^1(M)$ is realized in $M^*$.\\
This is easy with $|S^1(M)| \leq \kappa^+$.\\
{\bf Construction 2:} For all $M \prec N$ from $K_\kappa$ and $M \prec M' \in K_{\kappa^+}$, there is some $N' = *(M, N, M') \in K_{\kappa^+}$ such that $N, M' \prec N'$ and all of $S^1(N)$ are realized in $N'$.\\
For each $p \in S^1(N)$, find some $a_p \in |\sea|$ that realizes it.  Then find some $N' \prec \sea$ that contains $\{ a_p : p \in S^1(N) \} \cup |M'| \cup |N|$ of size $\kappa^+$.  This is possible since $|S^1(N)| \leq \kappa^+$.\\

{\bf Construction 3:} For all $M \in K_{\kappa^+}$ there is some $M^+ \in K_{\kappa^+}$ such that $M \prec M^+$ and, if $M_0 \prec M$ of size $\kappa$, then all of $S^1(M_0)$ are realized in $M^+$.\\
Find a resolution $\seq{M_i : i < \kappa^+}$ of $M$.  Set $N_0$ = $(M_0)^*$, $N_{i+1} = *(M_i, M_{i+1}, N_i)$, and take unions at limits.  Then $M^+ = \cup_{i < \kappa^+} N_i$ works.\\
{\bf Construction 4:} For all $M \in K_{\kappa^+}$, there is some $M^\# \in K_{\kappa^+}$ such that $M \prec M^\#$ and $M^\#$ is saturated.\\
Let $M \in K_\kappa$.  Set $M_0 = M$, $M_{i+1} = (M_i)^+$, and take unions at limits.  Then $M^\# = M_{\kappa^+}$ is saturated.\\
Then, to prove the proposition, let $M \in K_\kappa$.  Since $K$ has no maximal model, it has an extension $M'$ in $K_{\kappa^+}$.  Then $(M')^\#$ is the desired saturated extension of $M$.
\hfill \dag

\section{Further Work} \label{futurework}

In working with Galois types, as we do here, the assumption of amalgamation simplifies the definitions and construction by making $E_{AT}$ already an equivalence relation in the definition of types; see Definition \ref{galoistype}.  Could we obtain some bound on the number of long types in the absence of amalgamation?

In this paper, we introduced the definitions of $S_{sep}$ and $S_{strsep}$.  There are several basic questions to explore.

First, is $S_{strsep}$ necessary?  That is, is there an AEC where the size of $S_{strsep}$ is well behaved, but $S_{sep}$ is chaotic as in Proposition \ref{badnaex}; or can the proof of Theorem \ref{strsepone} be improved to provide the same bound on $seq\mathfrak{tb}$?  By Proposition \ref{dapsep}, any example of this chaotic behavior must have amalgamation but not disjoint amalgamation.  Perhaps one of the examples from \cite{bks} can be refined for this purpose.

In the original definitions of separative and strongly separative types, any ordered set $I$ was allowed as an index set and the separation properties were required to hold for all subsets instead of just an initial segment.  We give the original definition here under the names \emph{unordered separative} and \emph{unordered strongly separative types}.
\begin{defin}
\begin{enumerate}
	\item[]

	\item For $M \in K$, define $gS^\alpha_{\te{usep}}(M) = \{ gtp(\seq{a_i : i \in I}/M, N) \in gS^I_{\te{na}}(M) : \te{ for all $I_0 \subset I$, there is}$ $\te{ some $M \prec N_{I_0} \prec N$ such that $a_i \in N_{I_0}$ iff $i \in I_0$} \}$.

	\item For $M \in K$, define $gS^I_{\te{ustrsep}}(M) = \{ p = gtp(\seq{a_i : i \in I}/M, N) \in gS^I_{\te{usep}}(M) :$ for all $I_0 \subset I$ and $M \prec N_{I_0} \prec N$, if with $a_i \in N_{I_0}$ iff $i \in I_0$, then for every $f:N_{I_0} \to N_1^+$ with $\|N_1^+\| = N_{I_0}$, there is some $g:N \to N_2^+$ such that $f \subset g$ and $g(\seq{a_i : i \in I - I_0}) \not \in N_1^+ \}$.

\end{enumerate}
\end{defin}
Are these definitions equivalent to the ones given in Section \ref{strongseptypes}, or is there some example of an AEC where the two notions are distinct?  This question will likely be clarified by a lower bound for $strsep\mathfrak{tb}$ or an example lacking amalgamation.

A more lofty goal would be to attempt to classify stable, DAP AECs by the possible values of $na\mathfrak{tb}_\lambda^1 = sep\mathfrak{tb}_\lambda^1$.  That is, we know that $na\mathfrak{tb}^1_\lambda$ is a cardinal between $1$ and $\lambda$ and that this controls the value of $sep\mathfrak{tb}^\kappa_\lambda$ for all $\kappa \leq \lambda$.  For each value in $[1, \lambda] \cap CARD$, does an AEC with DAP and exactly that many non-algebraic types of length one exist?  Of particular interest is the discussion after Theorem \ref{strsepone}.  The conclusion of only one separative type (or strongly separative type, if we wish to drop the assumption of DAP) of any length over a model seems to be a very powerful hypothesis.

Looking back to the first order case, it would be interesting to find a syntactic characterization of separative types, keeping in mind that separative and strongly separative types are equivalent in this context.

\end{document}